\documentclass[a4paper]{article}
\usepackage{amssymb}
\usepackage{amsfonts}
\usepackage{amsmath}

\input{tcilatex}

\begin{document}

\title{Rational homotopy theory of flag manifolds}
\author{Haibao Duan\thanks{%
This work is supported by the National Science Foundation of China under
Grant no. 12331003.}}
\maketitle

\begin{abstract}
Let $G/P$ be a flag manifold with $G$ simple. We determine the rational
homotopy groups $\pi _{\ast }(G/P)\otimes \mathbb{Q}$ and describe the
rational cohomology ring $H^{\ast }(G/P;\mathbb{Q})$ solely in terms of the
rational homotopy types of $G$ and $P$.

The proof is based on the restrictions of basic Weyl invariants of $G$ to
suitable simple factors of $P$, thereby resolving a key nonvanishing issue
left open by earlier approaches of Meier, Shiga and Tezuka, Kotschick, and
Terzi\'{c}.
\end{abstract}

\section{Main results}

Let $G$ be a compact, connected Lie group with a maximal torus $T\subset G$.
The dimension $n:=\dim T$ is an invariant of $G$, called its \textsl{rank}.
Serre \cite{Se} showed that there exists a sequence of $n$ integers

\begin{center}
$\mathcal{F}_{G}=\left\{ r_{1},\cdots ,r_{n}\right\} $, $1\leq r_{1}\leq
\cdots \leq r_{n}$,
\end{center}

\noindent such that $G$ is rationally homotopy equivalent to a product of $n$
odd-dimensional spheres:

\begin{enumerate}
\item[(1.1)] $G\simeq _{\mathbb{Q}}S^{2r_{1}-1}\times \cdots \times
S^{2r_{n}-1}$.
\end{enumerate}

\noindent We call $\mathcal{F}_{G}$ the\textsl{\ rational homotopy type} of $%
G$.

Identify the Lie algebra $L(G)$ of $G$ with the tangent space of $G$ at the
identity element $e\in G$. Equip $L(G)$ with an inner product invariant
under the adjoint action of $G$. Let

\begin{center}
$\exp :L(G)\rightarrow G$
\end{center}

\noindent denote the exponential map. For a nonzero element $\alpha \in
L(G)\backslash \left\{ 0\right\} $ let $P_{\alpha }$ be the centralizer of
the one-parameter subgroup

\begin{center}
$\left\{ \exp (t\alpha )\in G\mid t\in \mathbb{R}\right\} $.
\end{center}

\noindent Then $P_{\alpha }$ is a parabolic subgroup, and the homogeneous
space $G/P_{\alpha }$ admits a canonical structure of a smooth complex
projective variety, called a \textsl{flag manifold} of $G$. Flag manifolds
constitute a fundamental class of compact homogeneous spaces. Their
cohomology has played a central role in the development of characteristic
classes, representation theory, and algebraic geometry, and has been
extensively studied through Schubert calculus; see, for example, \cite%
{BGP,DZ3}. The purpose of this paper is to establish a unified formula for
the rational homotopy groups of flag manifolds.

We first assume that $P$ is a parabolic subgroup of a simply connected
simple Lie group $G$. Write

\begin{center}
$\mathcal{F}_{G}\cap \mathcal{F}_{P}=\left\{ r_{1},\cdots ,r_{k}\right\} $, $%
r_{1}\leq \cdots \leq r_{k}$.
\end{center}

\noindent After reordering the remaining elements, we may write

\begin{enumerate}
\item[(1.2)] $\mathcal{F}_{G}=\left\{ r_{1},\cdots ,r_{k}\right\} \sqcup
\{j_{1},\cdots ,j_{n-k}\}$, $j_{1}\leq \cdots \leq j_{n-k}$,

$\mathcal{F}_{P}=\left\{ r_{1},\cdots ,r_{k}\right\} \sqcup \{i_{1},\cdots
,i_{n-k}\}$, $i_{1}\leq \cdots \leq i_{n-k}$.
\end{enumerate}

\noindent Here, the intersection is understood with multiplicities.

The following theorem establishes a fundamental principle in the rational
homotopy theory of flag manifolds: the rational homotopy groups of $G/P$ are
completely determined by the rational homotopy types $\mathcal{F}_{G}$ and $%
\mathcal{F}_{P}$ of $G$ and $P$, respectively.

\bigskip

\noindent \textbf{Theorem 1.1. }The rational homotopy groups of $G/P$ are
given by

\begin{enumerate}
\item[(1.3)] $\pi _{\ast }(G/P)\otimes \mathbb{Q}=\overset{n-k}{\underset{s=1%
}{\oplus }}(\mathbb{Q}[2i_{s}]\oplus \mathbb{Q}[2j_{s}-1])$,
\end{enumerate}

\noindent where $\mathbb{Q}[s]\subset \pi _{s}(G/P)\otimes \mathbb{Q}$
denotes a one-dimensional direct summand isomorphic to $\mathbb{Q}$.

\bigskip

\noindent \textbf{Example 1.2. }The conclusion of Theorem 1.1 does not
extend to arbitrary simply connected homogeneous spaces $G/K$, even when
both $G$ and $K$ are compact, connected, simply connected simple Lie groups.
To see this, consider the inclusions

\begin{center}
$K=SU(3)\subset SO(8)\subset SU(8)=G$,
\end{center}

\noindent where $SU(3)$ is embedded in $SO(8)$ via its adjoint
representation. Since

\begin{center}
$\mathcal{F}_{G}=\{2,3,4,5,6,7,8\},$ $\mathcal{F}_{K}=\{2,3\},$
\end{center}

\noindent we have

\begin{center}
$\mathcal{F}_{G}\cap \mathcal{F}_{K}=\mathcal{F}_{K}=\{2,3\}$.
\end{center}

\noindent Nevertheless, the inclusion $K\rightarrow G$ factors through $%
SO(8) $, and

\begin{center}
$\pi _{5}(SO(8))\otimes \mathbb{Q}=0$.
\end{center}

\noindent Consequently, the induced homomorphism

\begin{center}
$\pi _{5}(K)\otimes \mathbb{Q\rightarrow }\pi _{5}(SU(8))\otimes \mathbb{Q}$
\end{center}

\noindent is trivial. It follows from the rational homotopy exact sequence
of the fibration

\begin{center}
$K\hookrightarrow G\rightarrow G/K$
\end{center}

\noindent that

\begin{center}
$\pi _{\ast }(G/K)\otimes \mathbb{Q}=\mathbb{Q}[6]\overset{8}{\underset{s=3}{%
\oplus }}\mathbb{Q}[2s-1]$.
\end{center}

\noindent In particular, the rational homotopy groups of $G/K$ cannot be
determined solely by the rational homotopy types $\mathcal{F}_{G}$ and $%
\mathcal{F}_{K}$ through a formula analogous to that in Theorem 1.1.$\square 
$

\bigskip

Although the rational homotopy groups of homogeneous spaces have been
studied extensively, notably in the works of Meier, Shiga and Tezuka,
Kotschick, and Terzi\'{c} \cite{Me,ST,KT,T}, we are not aware of any formula
expressing the rational homotopy groups of a flag manifold $G/P$ solely in
terms of the rational homotopy types of $G$ and $P$, as provided by Theorem
1.1. See also Remark 3.5 for a comparison with these earlier approaches. 

By a theorem of Bott--Samelson \cite{BSa}, the flag manifold $G/P$ is simply
connected, and its integral cohomology is torsion free and concentrated in
even degrees. Together with (1.3), we find that $G/P$ is \textsl{rationally
elliptic} in the sense of \cite[p.434]{FHT}; that is

\begin{center}
$\dim H^{\ast }(G/P;\mathbb{Q})<\infty $, $\dim \pi _{\ast }(G/P)\otimes 
\mathbb{Q}<\infty $.

$\dim \pi _{odd}(G/P)\otimes \mathbb{Q}=\dim \pi _{even}(G/P)\otimes \mathbb{%
Q}=n-k$.
\end{center}

\noindent It follows from the F\'{e}lix--Halperin--Thomas theorem \cite[%
Proposition 32.10]{FHT} that the rational cohomology ring of $G/P$ admits
the following description in terms of the rational homotopy types $\mathcal{F%
}_{G}$ and $\mathcal{F}_{P}$.

\bigskip

\noindent \textbf{Corollary 1.3. }Let $G$ be a simply connected simple Lie
group with a parabolic subgroup $P\subset G$, and suppose that the rational
homotopy types of $G$ and $P$ are given by (1.2). Then

\begin{center}
$H^{\ast }(G/P;\mathbb{Q})\cong \frac{\mathbb{Q[}x_{1},\cdots ,x_{n-k}%
\mathbb{]}}{\left\langle f_{1},\cdots ,f_{n-k}\right\rangle }$, $\deg
x_{s}=2i_{s}$,
\end{center}

\noindent where $f_{1},\cdots ,f_{n-k}$ form a regular sequence of
decomposable polynomials in $\mathbb{Q[}x_{1},\cdots ,x_{n-k}\mathbb{]}$
satisfying

\begin{center}
$\deg f_{s}=2j_{s},1\leq s\leq n-k$.
\end{center}

Finally, since $G$ is simple, (1.2) implies that the degrees of the
generators and those of the defining relations are disjoint; namely,

\begin{enumerate}
\item[(1.4)] $\{i_{1},\cdots ,i_{n-k}\}\cap \{j_{1},\cdots
,j_{n-k}\}=\emptyset $.
\end{enumerate}

\bigskip

\noindent \textbf{Remark 1.4. }The rational homotopy types of simple Lie
groups are well known and are listed in Table 1 in Section 2. Moreover,
Theorem 2.4 and Corollary 2.2 provide an explicit procedure for determining $%
\mathcal{F}_{P}$ for each parabolic subgroup $P$; see Example 2.5.
Consequently, Theorem 1.1 and Corollary 1.3 can be applied directly to
determine the rational homotopy groups and the rational cohomology ring of
the flag manifold $G/T$.

More generally, every flag manifold $G/P$ admits a product decomposition, up
to diffeomorphism,

\begin{enumerate}
\item[(1.5)] $G/P\cong G_{1}/P_{1}\times \cdots \times G_{m}/P_{m}$,
\end{enumerate}

\noindent where each $G_{i}$ is simply connected and simple, and $P_{i}$ is
a parabolic subgroup of $G_{i}$. Thus, Theorem 1.1 and Corollary 1.3 can be
applied factorwise to describe the rational homotopy groups and rational
cohomology ring of an arbitrary flag manifold.

We also note that the disjointness property (1.4) is specific to the case
where $G$ is simple. For a nonsimple $G$, this property may fail, since the
product decomposition (1.5) can give rise to common degrees among the
algebra generators and defining relations of the rational cohomology ring.$%
\square $

\bigskip

\noindent \textbf{Remark 1.5. }According to A. Weil \cite[p.331]{W},
classical Schubert calculus amounts to determining the integral cohomology
ring $H^{\ast }(G/P)$ of a flag manifold $G/P$. Using a multiplicative rule 
\cite{D} for Schubert classes, Duan and Zhao developed an algorithm for
computing $H^{\ast }(G/P)$ in terms of explicit generators and defining
relations expressed in Schubert classes; see, for example, \cite{DZ1,DZ2}.
In this approach, Corollary 1.3 provides a useful algebraic guide to the
computation by determining the degrees of the generators and defining
relations in advance.$\square $

\bigskip

\noindent \textbf{Remark 1.6.} For two topological spaces $X$ and $Y$, let $%
[X,Y]$ denote the set of homotopy classes of maps from $X$ to $Y$, and let $%
\mathrm{Hom}(H^{\ast }(Y),H^{\ast }(X))$ denote the set of degree-preserving
ring homomorphisms from $H^{\ast }(Y)$ to $H^{\ast }(X)$. Assigning to a map 
$f:X\rightarrow Y$ its induced homomorphism $f^{\ast }$ defines a natural map

\begin{center}
$R:[X,Y]\rightarrow \mathrm{Hom}(H^{\ast }(Y),H^{\ast }(X))$, $R[f]=f^{\ast
} $.
\end{center}

\noindent When both $X$ and $Y$ are flag manifolds, the map $R$ has finite
kernel and cokernel; see \cite{HM,Sh,P}. More precisely, let $X_{0}$ and $%
Y_{0}$ denote the rationalizations of $X$ and $Y$, respectively. Since flag
manifolds are simply connected and formal, Sullivan's theory \cite{S} yields
a natural bijection \cite{GH,P}

\begin{center}
$[X_{0},Y_{0}]\cong \mathrm{Hom}(H^{\ast }(Y;\mathbb{Q}),H^{\ast }(X;\mathbb{%
Q}))$.
\end{center}

\noindent Thus, for flag manifolds, the rational cohomology ring determines
not only the rational homotopy type of each manifold but also the homotopy
classes of maps between their rationalizations. In particular, the
presentation of $H^{\ast }(G/P;\mathbb{Q})$ given in Corollary 1.3 provides
an effective algebraic framework for studying maps between rationalized flag
manifolds. For applications of this approach to the classification of
self-maps of certain flag manifolds, see, for example, \cite{DHZ,P}.$\square 
$

\bigskip

The remainder of the paper is organized as follows. Section 2 recalls the
classification of parabolic subgroups $P$ of a simple Lie group $G$,
computes their rational homotopy types $\mathcal{F}_{P}$, and identifies the
possible simple factors of their semisimple parts. Section 3 studies the
restrictions of basic Weyl invariants of $G$ to these simple factors,
thereby determining the induced map of the inclusion $i:P\rightarrow G$ on
rational homotopy groups. Building on these results, we complete the proof
of Theorem 1.1 in Section 5.

\section{Classification of parabolic subgroups of a simple Lie group}

The simply connected simple Lie groups fall into four infinite families of
classical types,

\begin{center}
$\left\{ A_{n},\ n\geq 1\right\} $, $\left\{ B_{n},\ n\geq 2\right\} $, $%
\left\{ C_{n},\ n\geq 3\right\} $, $\left\{ D_{n},\ n\geq 4\right\} $,
\end{center}

\noindent and five exceptional types:

\begin{center}
$G_{2},F_{4},E_{6},E_{7},E_{8}$.
\end{center}

\noindent Their rational homotopy types $\mathcal{F}_{G}$ and centers $%
\mathcal{Z}(G)$ are listed in Table 1; see \cite[Table 1]{K} and \cite[p.57]%
{Hu}.

\begin{center}
\begin{tabular}{|l||l|l|}
\hline
$G$ & $\mathcal{F}_{G}$ & $\mathcal{Z}(G)$ \\ \hline\hline
$A_{n}$ & $\left\{ 2,3,4,\cdots ,n+1\right\} $ & $\mathbb{Z}_{n+1}$ \\ \hline
$B_{n}$ & $\left\{ 2,4,6,\cdots ,2n\right\} $ & $\mathbb{Z}_{2}$ \\ \hline
$C_{n}$ & $\left\{ 2,4,6,\cdots ,2n\right\} $ & $\mathbb{Z}_{2}$ \\ \hline
$D_{n}$ & $\left\{ 2,4,\cdots ,2n-2,n\right\} $ & $\mathbb{Z}_{4}${\small ,\ 
}${\small n}$ odd; $\mathbb{Z}_{2}{\small \oplus }\mathbb{Z}_{2}${\small , }$%
{\small n}$ even \\ \hline
$G_{2}$ & $\left\{ 2,6\right\} $ & ${\small \{e\}}$ \\ \hline
$F_{4}$ & $\left\{ 2,6,8,12\right\} $ & ${\small \{e\}}$ \\ \hline
$E_{6}$ & $\left\{ 2,5,6,8,9,12\right\} $ & $\mathbb{Z}_{3}$ \\ \hline
$E_{7}$ & $\left\{ 2,6,8,10,12,14,18\right\} $ & $\mathbb{Z}_{2}$ \\ \hline
$E_{8}$ & $\left\{ 2,8,12,14,18,20,24,30\right\} $ & ${\small \{e\}}$ \\ 
\hline
\end{tabular}

{\small Table 1. The rational homotopy types and centers of simply-connected
simple Lie groups}
\end{center}

\noindent In general, every compact connected Lie group $P$ admits a
presentation

\begin{enumerate}
\item[(2.1)] $P=(G_{1}\times \cdots \times G_{k}\times T^{m})/F$
\end{enumerate}

\noindent where each factor $G_{t}$ is a simply connected simple Lie groups
listed in Table 1, $T^{m}$ is a torus of rank $m$, and $F$ is a finite
subgroup of the product

\begin{center}
$\mathcal{Z}(G_{1})\times \cdots \times \mathcal{Z}(G_{k})\times T^{m}$.
\end{center}

\noindent Motivated by this presentation, we introduce the following
terminology.

\bigskip

\noindent \textbf{Definition 2.1.} The numerator $G_{1}\times \cdots \times
G_{k}\times T^{m}$ in (2.1) is called the \textsl{local type} of $P$. The
product

\begin{center}
$P^{s}:=G_{1}\times \cdots \times G_{k}$
\end{center}

\noindent is called the \textsl{semisimple part} of $P$. Each factor $G_{t}$
is called a \textsl{simple factor} of the semisimple part of $P$.\noindent $%
\square $

\bigskip

This follows from the fact that the natural projection in (2.1) is a finite
covering map.

\bigskip

\noindent \textbf{Corollary 2.2. }The rational homotopy type of the group $P$
is determined by its local type and is given by

\begin{center}
$\mathcal{F}_{P}=\left\{ 1,\cdots ,1\right\} $($m$ times) $\sqcup \mathcal{F}%
_{G_{1}}\sqcup \cdots \sqcup \mathcal{F}_{G_{k}}$.$\square $
\end{center}

\bigskip

Thus, in order to apply Theorem 1.1 effectively to the computation of the
rational homotopy groups $\pi _{\ast }(G/P)\otimes \mathbb{Q}$, it suffices
to classify the parabolic subgroups $P$ of a simple Lie group $G$ in terms
of their local types.

Fix a maximal torus $T\subset G$. Using the invariant inner product to
identify the dual space $L(T)^{\ast }$ with $L(T)$, let

\begin{center}
$\Omega =\{\omega _{{1}},\ldots ,\omega _{{n}}\}\subset L(T)$
\end{center}

\noindent be a set of \textsl{fundamental dominant weights }\cite[p.67]{Hu}.
Geometrically, the closed convex cone generated by $\omega _{{1}},\ldots
,\omega _{{n}}$

\begin{center}
$\mathcal{K}=\left\{ \lambda _{1}\omega _{{1}}+\cdots +\lambda _{n}\omega _{{%
n}}\in L(T)\mid \lambda _{1},\cdots ,\lambda _{n}\geq 0\right\} $,
\end{center}

\noindent can be identified with a closed Weyl chamber of $G$, whose extreme
rays are spanned by $\omega _{{1}},\ldots ,\omega _{{n}}$.

For each subset $I\subseteq \{1,\cdots ,n\}$, let $P_{I}$ denote the
centralizer of the one-parameter subgroup

\begin{center}
$\alpha _{I}:\mathbb{R}\rightarrow G$, $\alpha _{I}(t)=\exp
(t\sum\limits_{i\in I}\omega _{i})$.
\end{center}

\noindent The following classification is known; see \cite[Lemma 2.2]{DZ2}.

\bigskip

\noindent \textbf{Lemma 2.3}.\textbf{\ }Every parabolic subgroup $P$ of $G$
is conjugate to $P_{I}$\ for some $I\subseteq \{1,\cdots ,n\}$.
Consequently, every flag manifold of $G$ is isomorphic to a homogeneous
space $G/P_{I}$.\noindent $\square $

\bigskip

Recall that simply connected semi-simple Lie groups $G$ are classified by
their Dynkin diagrams $\Gamma _{G}$ (see \cite[p.58]{Hu}). Throughout, for a
simply connected simple Lie group $G$, we identify the fundamental dominant
weights in $%
\Omega
$ with the vertices of $\Gamma _{G}$, and label these vertices accordingly,
following the convention in \cite[p.58]{Hu}.

Note that a parabolic subgroup $P$ of a simple Lie group $G$ is a connected
subgroup of maximal rank. The following description of the local types of
parabolic subgroups is a consequence of the Borel--Siebenthal classification 
\cite{BS} of connected subgroups of maximal rank in compact simple Lie
groups.

\bigskip

\noindent \textbf{Theorem 2.4. }Let $G$ be a simply connected simple Lie
group of rank $n$,\textbf{\ }with Dynkin diagram $\Gamma _{G}$.\textbf{\ }%
For a subset $I\subseteq \{1,\cdots ,n\}$, let $P_{I}$ be the corresponding
parabolic subgroup, and denote its local type by $P_{I}^{s}\times T^{m}$.
Then:

(i) $m=\left\vert I\right\vert $, the cardinality of $I$;

(ii) The Dynkin diagram of the semisimple part $P_{I}^{s}$ is obtained from $%
\Gamma _{G}$ by deleting the vertices indexed by $I$, together with all
edges incident to these vertices.$\square $

\bigskip

\noindent \textbf{Example 2.5. }Theorem 2.4, together with Corollary 2.2%
\textbf{, }provides an effective procedure for computing the rational
homotopy types\textbf{\ }$\mathcal{F}_{P_{I}}$ of the parabolic subgroups $%
P_{I}$, in term of the subset $I\subseteq \{1,\cdots ,n\}$. These
computations, in turn, allow us to apply formula (1.3) to determine the
rational homotopy groups of the corresponding flag manifolds $G/P_{I}$.

As an illustration, consider the Dynkin diagram of the exceptional Lie group 
$G=E_{6}$:

\begin{center}
$%
\begin{array}{ccccccccc}
&  &  &  & \omega _{2} &  &  &  &  \\ 
&  &  &  & \mid  &  &  &  &  \\ 
\omega _{1} & - & \omega _{3} & - & \omega _{4} & - & \omega _{5} & - & 
\omega _{6}%
\end{array}%
$.
\end{center}

\noindent For each vertex $\omega _{i}$, $1\leq i\leq 6$, take $I=\{i\}$. By
Theorem 2.4, the semisimple parts $P_{I}^{s}$ of the corresponding parabolic
subgroups are listed in the second column of Table 2. Their rational
homotopy types $\mathcal{F}_{P_{I}}$, obtained from Corollary 2.2, are given
in the third column.

\begin{center}
\begin{tabular}{l||l|l|}
\hline
$I$ & $P_{I}^{s}$ & $\mathcal{F}_{P_{I}}$ \\ \hline\hline
$\left\{ 1\right\} $ & $D_{5}$ & $\{1\}\sqcup \left\{ 2,4,6,8,5\right\} $ \\ 
\hline
$\left\{ 2\right\} $ & $A_{5}$ & $\{1\}\sqcup \left\{ 2,3,4,5,6\right\} $ \\ 
\hline
$\left\{ 3\right\} $ & $A_{1}\times A_{4}$ & $\{1\}\sqcup \left\{ 2\right\}
\sqcup \left\{ 2,3,4,5\right\} $ \\ \hline
$\left\{ 4\right\} $ & $A_{2}\times A_{1}\times A_{2}$ & $\{1\}\sqcup
\left\{ 2,3\right\} \sqcup \left\{ 2\right\} \sqcup \left\{ 2,3\right\} $ \\ 
\hline
$\left\{ 5\right\} $ & $A_{1}\times A_{4}$ & $\{1\}\sqcup \left\{ 2\right\}
\sqcup \left\{ 2,3,4,5\right\} $ \\ \hline
$\left\{ 6\right\} $ & $D_{5}$ & $\{1\}\sqcup \left\{ 2,4,6,8,5\right\} $ \\ 
\hline
\end{tabular}

{\small Table 2. The rational homotopy types of the parabolic subgroups }$%
{\small P}_{\left\{ i\right\} }$ {\small of }${\small E}_{6}$.$\square $
\end{center}

According to Lemma 2.3, a simple Lie group $G$ of rank $n$ has $2^{n}-1$
proper standard parabolic subgroups, in one-to-one correspondence with the
subsets of vertices of its Dynkin diagram $\Gamma _{G}$. Thus, a direct
proof of Theorem 1.1 might seem to require a large number of individual
cases. The following consequence of Theorem 2.4(ii), however, considerably
reduces the number of cases and provides a systematic framework for the
proof.

\bigskip

\noindent \textbf{Corollary 2.6.} Let $G$ be a simple Lie group. The
possible isomorphism types of the simple factors occurring in the semisimple
parts of its proper standard parabolic subgroups are precisely those listed
in Table 3.

\begin{center}
\begin{tabular}{l||l|}
\hline
$G$ & Possible simple factors $L$ \\ \hline\hline
$A_{n}$ & $A_{m}$, $1\leq m\leq n-1$ \\ \hline
$B_{n}$ & $B_{m}$, $2\leq m\leq n-1$; $A_{m}$, $1\leq m\leq n-1$ \\ \hline
$C_{n}$ & $C_{m}$, $2\leq m\leq n-1$; $A_{m}$, $1\leq m\leq n-1$ \\ \hline
$D_{n}$ & $D_{m}$, $4\leq m\leq n-1$; $A_{m}$, $1\leq m\leq n-1$ \\ \hline
$G_{2}$ & $A_{1}$ \\ \hline
$F_{4}$ & $A_{m}$, $1\leq m\leq 2$; $B_{3},C_{3}$ \\ \hline
$E_{6}$ & $A_{m}$, $1\leq m\leq 5$; $D_{k}$, $4\leq k\leq 5$ \\ \hline
$E_{7}$ & $E_{6}$; $A_{m}$, $1\leq m\leq 6$; $D_{k}$, $4\leq k\leq 6$ \\ 
\hline
$E_{8}$ & $E_{6}$; $E_{7}$; $A_{m}$, $1\leq m\leq 7$; $D_{k}$, $4\leq k\leq
7 $ \\ \hline
\end{tabular}%
.$\square $

{\small Table 3. Possible simple factors in the semisimple parts of
parabolic subgroups of simple Lie groups.}$\square $
\end{center}

\noindent \textbf{Remark 2.7.} In \cite{BS}, Borel and Siebenthal classified
maximal connected subgroups of maximal rank in compact simple Lie groups, up
to local type. Their classification provides a description of the local
types relevant to the parabolic subgroups considered here. In \cite[Theorem
4.4]{DL}, Duan and Liu refined this classification by determining the
corresponding isomorphism types, in the sense of the presentation (2.1).$%
\square $

\section{Restrictions of basic Weyl invariants}

Let $A=\underset{r\geq 0}{\oplus }A^{r}$ be a graded commutative algebra
over the rational field $\mathbb{Q}$, with $A^{0}=\mathbb{Q}$. Let

\begin{center}
$A^{+}=\underset{r\geq 1}{\oplus }A^{r}$
\end{center}

\noindent be the augmentation ideal. The \textsl{ideal of decomposable
elements} is the square $(A^{+})^{2}$. The quotient

\begin{center}
$Q(A)=A^{+}/(A^{+})^{2}$
\end{center}

\noindent is called \textsl{the space of indecomposables} of $A$. It is a
graded vector space, whose degree $r$ component is

\begin{center}
$Q^{r}(A)=A^{r}/\left( \dsum\limits_{i+j=r,\text{ }i,j>0}A^{i}A^{j}\right) $.
\end{center}

\noindent We denote by

\begin{center}
$q_{A}:A^{r}\rightarrow Q^{r}(A)$
\end{center}

\noindent the canonical quotient map.

\bigskip

\noindent \textbf{Example 3.1.} Let $B_{P}$ denote the classifying space of
a compact connected Lie group $P$. When $A=H^{\ast }(B_{P};\mathbb{Q})$ we
use the abbreviations

\begin{center}
$q_{P}:=q_{A}$, $\ Q^{r}(P):=Q^{r}(A)$.
\end{center}

\noindent Recall that if the rational homotopy type of $P$ is

\begin{center}
$\mathcal{F}_{P}=\left\{ r_{1},\cdots ,r_{n}\right\} $, $n=\mathrm{\limfunc{%
rank}}P$,
\end{center}

\noindent then $A$ is a polynomial algebra on $n$ homogeneous generators:

\begin{center}
$A=\mathbb{Q}[I_{r_{1}},\cdots ,I_{r_{n}}]$, $\deg I_{r_{k}}=2r_{k}$, $1\leq
k\leq n$.
\end{center}

\noindent Moreover, by Borel \cite{B}, these generators can be chosen to
form a basic set of homogeneous Weyl invariants of $P$. In particular,

\begin{center}
$\left\{ q_{P}(I_{r_{1}}),\cdots ,q_{P}(I_{r_{n}})\right\} $
\end{center}

\noindent is a homogeneous basis of\textsl{\ }$Q(P)$.

For a subgroup $L\subset P$ and an element $I\in H^{\ast }(B_{P};\mathbb{Q})$%
, we denote by $I\mid _{L}$ its image under the restriction homomorphism

\begin{center}
$H^{\ast }(B_{P};\mathbb{Q})\rightarrow H^{\ast }(B_{L};\mathbb{Q})$
\end{center}

\noindent induced by the inclusion $L\subset P$.$\square $

\bigskip

Let $G$ be a compact, connected, simple Lie group, and let $P\subset $ $G$
be a parabolic subgroup. As in Section 1, let $\mathcal{F}_{G}$ and $%
\mathcal{F}_{P}$ denote the rational homotopy types of $G$ and $P$,
respectively. Fix a basic set of homogeneous Weyl invariants

\begin{center}
$\left\{ I_{r_{1}},\cdots ,I_{r_{n}}\right\} $
\end{center}

\noindent of $G$, where $\mathcal{F}_{G}=\left\{ r_{1},\cdots ,r_{n}\right\} 
$. Recall that $r_{1}=2$ for every simple Lie group $G$.

\bigskip

\noindent \textbf{Lemma 3.2}. For every $r\in \mathcal{F}_{G}\cap \mathcal{F}%
_{P}$, there exists a basic Weyl invariant $I\in H^{2r}(B_{G};\mathbb{Q})$
such that

\begin{enumerate}
\item[(3.1)] $q_{P}(I\mid _{P})\neq 0$ in $Q^{2r}(P)$.
\end{enumerate}

\bigskip

\noindent \textbf{Observation 3.3}. Although $G$ has $2^{n}-1$ standard
parabolic subgroups, the following observations, together with the
classification in Corollary 2.6, substantially reduce the number of cases
that need to be considered in proving Lemma 3.2.

(i) Suppose that $r\in \mathcal{F}_{G}\cap \mathcal{F}_{P}$. By Theorem 2.4
and Corollary 2.2, the semisimple part of the local type of $P$ has a simple
factor $L$ such that $r\in \mathcal{F}_{L}$. It therefore suffices to find a
basic Weyl invariant $I\in H^{2r}(B_{G};\mathbb{Q})$ satisfying

\begin{center}
$q_{L}(I\mid _{L})\neq 0$ in $Q^{2r}(L)$.
\end{center}

(ii) Since $r_{1}=2$ for every simple Lie group $G$, the invariant $%
I_{r_{1}} $ is the quadratic invariant induced by the Killing form. Its
restriction to every simple factor of the local type of $P$ is nonzero in $%
Q^{4}(L)$. Thus, in proving Lemma 3.2, we may assume that $r>2$.

(iii) If $L_{1}\subset L_{2}$ are simple subgroups under consideration, let

\begin{center}
$q_{L_{2},L_{1}}:Q^{\ast }(L_{2})\rightarrow Q^{\ast }(L_{1})$
\end{center}

\noindent denote the homomorphism induced by the restriction

\begin{center}
$H^{\ast }(B_{L_{2}};\mathbb{Q})\rightarrow H^{\ast }(B_{L_{1}};\mathbb{Q})$.
\end{center}

\noindent Then

\begin{center}
$q_{L_{1}}(I\mid _{L_{1}})=q_{L_{2},L_{1}}(q_{L_{2}}(I\mid _{L_{2}}))$.
\end{center}

\noindent In particular, if $q_{L_{1}}(I\mid _{L_{1}})\neq 0$, then $%
q_{L_{2}}(I\mid _{L_{2}})\neq 0$.$\square $

\bigskip

\noindent \textbf{Proof of Lemma 3.2. }Suppose that $r\in \mathcal{F}%
_{G}\cap \mathcal{F}_{P}$. In view of Observation 3.3, we distinguish two
cases according to whether $G$ is of classical or exceptional type.

For the classical simple groups, the relevant restriction formulas are
summarized in Table 4. The basic Weyl invariants

\begin{center}
$I_{r}\in H^{2r}(B_{G};\mathbb{Q})$ and $p_{r}\in H^{2r}(B_{L};\mathbb{Q})$,
\end{center}

\noindent appearing in the second and fifth columns, respectively, are
chosen to be the standard power-sum generators of the corresponding
invariant rings. For type $D_{n}$, we also need the Euler invariant

\begin{center}
$e_{n}=x_{1}\cdots x_{n}\in $ $H^{2n}(B_{D_{n}};\mathbb{Q})$.
\end{center}

\noindent When $n$ is even, $e_{n}$ is an additional basic invariant of
degree $2n$, alongside the power-sum invariant $I_{n}$; when $n$ is odd, $%
e_{n}$ is the basic invariant of degree $2n$.

The last column of Table 4 records the corresponding nonzero indecomposable
classes

\begin{center}
$q_{L}(I_{r}^{\prime }\mid _{L})\in Q^{2r}(L)$,
\end{center}

\noindent where

\begin{center}
$I_{r}^{\prime }=\left\{ 
\begin{tabular}{l}
$e_{n}\text{, if }G=D_{n}\text{, }L=A_{n-1}\text{, }r=n\text{, }n\text{ odd,}
$ \\ 
$I_{r}\text{, otherwise}$.%
\end{tabular}%
\right. $
\end{center}

\noindent In this exceptional case, the restriction satisfies

\begin{center}
$q_{A_{n-1}}(e_{n}\mid _{A_{n-1}})=\frac{1}{n}[p_{n}]\neq 0$.
\end{center}

\noindent In view of the classification of the possible simple factors $L$
of $P$ given in Corollary 2.6, Table 4 verifies Lemma 3.2 for all classical
simple groups.

\begin{center}
\begin{tabular}{l||l|l|l|l}
\hline
$G$ & Basic Weyl invariants & $L$ & $r\in \mathcal{F}_{G}\cap \mathcal{F}%
_{L} $ & $q_{L}(I_{r}^{\prime }\mid _{L})\in Q^{2r}(L)$ \\ \hline\hline
$A_{n}$ & $I_{2},I_{3},\cdots ,I_{n+1}$ & $A_{m}$ & $2\leq r\leq m+1$ & $%
\left[ p_{r}\right] $ \\ \hline
$B_{n}$ & $I_{2},I_{4},\cdots ,I_{2n}$ & $B_{m}$ & $r=2j\text{, }1\leq j\leq
m$ & $\left[ p_{r}\right] $ \\ \hline
$B_{n}$ & $I_{2},I_{4},\cdots ,I_{2n}$ & $A_{m}$ & $r=2j$, $1\leq j\leq
\left\lfloor \frac{m+1}{2}\right\rfloor $ & $\left[ p_{r}\right] $ \\ \hline
$C_{n}$ & $I_{2},I_{4},\cdots ,I_{2n}$ & $C_{m}$ & $r=2j\text{, }1\leq j\leq
m$ & $\left[ p_{r}\right] $ \\ \hline
$C_{n}$ & $I_{2},I_{4},\cdots ,I_{2n}$ & $A_{m}$ & $r=2j\text{, }1\leq j\leq
\left\lfloor \frac{m+1}{2}\right\rfloor $ & $\left[ p_{r}\right] $ \\ \hline
$D_{n}$ & $I_{2},I_{4},\cdots ,I_{2n-2},e_{n}$ & $D_{m}$ & $r=2j\text{, }%
1\leq j\leq n-1$ & $\left[ p_{r}\right] $ \\ \hline
$D_{n}$ & $I_{2},I_{4},\cdots ,I_{2n-2},e_{n}$ & $A_{m}$ & $r=2j\text{, }%
2\leq r\leq n-1$ & $\left[ p_{r}\right] $ \\ \hline
$D_{n}$ & $e_{n}$ & $A_{n-1}$ & $r=n$, $n$ odd & $\frac{1}{n}[p_{n}]$ \\ 
\hline
\end{tabular}%
.

{\small Table 4. Restrictions of basic Weyl invariants for classical simple
groups.}
\end{center}

Turning to the exceptional groups, we use the basic Weyl invariants

\begin{center}
$I_{r_{1}},\cdots ,I_{r_{n}}$, $\mathcal{F}_{G}=\left\{ r_{1},\cdots
,r_{n}\right\} $,
\end{center}

\noindent given by Mehta \cite{M}. By Corollary 2.6, the simple factors that
can occur in the semisimple parts of the local types of parabolic subgroups
are explicitly listed. Restricting the basic Weyl invariants of $G$ to these
factors and passing to the corresponding indecomposable quotients yields the
formulas collected in Tables 5--8, which verify Lemma 3.2 for $%
G=F_{4},E_{6},E_{7}$ and $E_{8}$, respectively.

\begin{center}
\begin{tabular}{l||l|l|l}
\hline
B$\text{asic Weyl invariants}$ & $L$ & $r\in \mathcal{F}_{G}\cap \mathcal{F}%
_{L}$ & $q_{L}(I_{r}\mid _{L})\in Q^{2r}(L)$ \\ \hline\hline
$I_{2},I_{6},I_{8},I_{12}$ & $B_{3}$ & $6$ & $\frac{9}{4}[p_{6}]$ \\ \hline
& $C_{3}{}$ & $6$ & $-48[p_{6}]$ \\ \hline
\end{tabular}

{\small Table 5. Restrictions of basic Weyl invariants for }${\small F}_{4}$

\bigskip

\begin{tabular}{l||l|l|l}
\hline
B$\text{asic Weyl invariants}$ & $L$ & $r\in \mathcal{F}_{G}\cap \mathcal{F}%
_{L}$ & $q_{L}(I_{r}\mid _{L})\in Q^{2r}(L)$ \\ \hline\hline
$I_{2},I_{5},I_{6},I_{8},I_{9},I_{12}$ & $A_{4}$ & $5$ & $-10[p_{5}]$ \\ 
\hline
& $A_{5}{}$ & $5,6$ & $-10[p_{5}],-24[p_{6}]$ \\ \hline
& $D_{4}$ & $6$ & $6[p_{6}]$ \\ \hline
& $D_{5}$ & $5,6,8$ & $60[e_{5}],6[p_{6}],-15[p_{8}]$ \\ \hline
\end{tabular}

{\small Table 6. Restrictions of basic Weyl invariants for }${\small E}_{6}$.
\end{center}

For $E_{7}$, we make the following remarks.

(i) When $L=E_{6}$, we denote Mehta's basic Weyl invariant of degree $2r$ by 
$J_{r}$.

(ii) The two types of $A_{5}$-factors must be distinguished, since they
belong to two distinct conjugacy classes of simple $A_{5}$-subgroups in $%
E_{7}$.

(iii) Since $\mathcal{F}_{D_{6}}=\{2,4,6,6,8,10\}$, the indecomposable
quotient $Q^{12}(D_{6})$ is two-dimensional. Accordingly, both $[p_{6}]$ and 
$[e_{6}]$ are needed to describe the restriction of $I_{6}$ to $D_{6}$.

\begin{center}
\begin{tabular}{l||l|l|l}
\hline
B$\text{asic Weyl invariants}$ & $L$ & $r\in \mathcal{F}_{G}\cap \mathcal{F}%
_{L}$ & $q_{L}(I_{r}\mid _{L})\in Q^{2r}(L)$ \\ \hline\hline
$I_{2},I_{6},I_{8},I_{10},I_{12},I_{14},I_{18}$ & $E_{6}$ & $6,8,12$ & $%
2[J_{6}],2[J_{8}],2[J_{12}]$ \\ \hline
& $D_{6}$ & $6,8,10$ & $\frac{9}{2}[p_{6}]-360[e_{6}],\frac{33}{8}[p_{8}],%
\frac{129}{32}[p_{10}]$ \\ \hline
& $D_{5}$ & $6,8$ & $\frac{9}{2}[p_{6}],\frac{33}{8}[p_{8}]$ \\ \hline
& $A_{6}$ & $6$ & $12[p_{6}]$ \\ \hline
& $A_{5}$ & $6$ & $\frac{9}{2}[p_{6}]$ \\ \hline
& $A_{5}^{\prime }$ & $6$ & $12[p_{6}]$ \\ \hline
& $D_{4}$ & $8$ & $\frac{33}{8}[p_{8}]$ \\ \hline
\end{tabular}

{\small Table 7. Restrictions of basic Weyl invariants for }${\small E}_{7}$.
\end{center}

For $E_{8}$, we denote Mehta's \cite{M} basic invariant of degree $2r$ for
the simple factors $E_{6},E_{7}\subset E_{8}$ by $K_{r}$ and $J_{r}$,
respectively.

\begin{center}
\begin{tabular}{l||l|l|l}
\hline
B$\text{asic Weyl invariants}$ & $L$ & $r\in \mathcal{F}_{G}\cap \mathcal{F}%
_{L}$ & $q_{L}(I_{r}\mid _{L})\in Q^{2r}(L)$ \\ \hline\hline
$I_{2},I_{8},I_{12},I_{14},I_{18},I_{20},I_{24},I_{30}$ & $E_{7}$ & $%
8,12,14,18$ & $6[J_{8}],-14[J_{12}],\frac{300}{29}[J_{14}],-\frac{54600}{1229%
}[J_{18}],$ \\ \hline
& $E_{6}$ & $8,12$ & $12[K_{8}],-28[K_{12}]$ \\ \hline
& $D_{7}$ & $8,12$ & $-180[p_{8}],-7560[p_{12}]$ \\ \hline
& $D_{6}$ & $8$ & $-180[p_{8}]$ \\ \hline
& $D_{5}$ & $8$ & $-180[p_{8}]$ \\ \hline
& $A_{7}$ & $8$ & $-1440[p_{8}]$ \\ \hline
\end{tabular}

{\small Table 8. Restrictions of basic Weyl invariants for }${\small E}_{8}$.
\end{center}

\noindent This completes the proof of Lemma 3.2. The nonzero classes
recorded in the last columns of Tables 4--8 establish the required
nonvanishing. The precise values of the coefficients are immaterial; only
their nonvanishing is used in the proof.$\square $

\bigskip

Let $G$ be a simply connected simple Lie group, and let $i:P\rightarrow G$
denote the inclusion of a parabolic subgroup. Lemma 3.2 can be reformulated
as a nonvanishing property of the induced homomorphism on rational homotopy
groups.

\bigskip

\noindent \textbf{Theorem 3.4. }For every $r\in \mathcal{F}_{G}\cap \mathcal{%
F}_{P}$, the induced map

\begin{center}
$i_{\ast }:\pi _{2r-1}(P)\otimes \mathbb{Q\rightarrow }\pi _{2r-1}(G)\otimes 
\mathbb{Q}$.
\end{center}

\noindent is nontrivial.\noindent

\bigskip

\noindent \textbf{Proof. }Recall that, for a compact connected Lie group $P$%
, rational homotopy theory provides a natural isomorphism

\begin{center}
$\mathcal{T}_{P}:Q^{s}(P)\cong \mathrm{Hom}(\pi _{s-1}(P),\mathbb{Q})$
\end{center}

\noindent where $Q^{s}(P)$ denotes the degree-$s$ component of the space of
indecomposables of $H^{\ast }(B_{P};\mathbb{Q})$; see, for example, F\'{e}%
lix-Halperin-Thomas \cite{FHT}. More precisely, for any homomorphism of
compact connected Lie groups $f:P\rightarrow G$, the induced homomorphism

\begin{center}
$Q^{s}(f):Q^{s}(G)\rightarrow Q^{s}(P)$
\end{center}

\noindent fits into the following commutative diagram

\begin{center}
\begin{tabular}{lll}
$Q^{s}(G)$ & $\overset{\cong }{\rightarrow }$ & $\mathrm{Hom}(\pi _{s-1}(G),%
\mathbb{Q})$ \\ 
$Q^{s}(f)\downarrow $ &  & $\downarrow \mathrm{Hom}(f_{\ast },\mathrm{id})$
\\ 
$Q^{s}(P)$ & $\overset{\cong }{\rightarrow }$ & $\mathrm{Hom}(\pi _{s-1}(P),%
\mathbb{Q})$%
\end{tabular}%
.
\end{center}

Now let $i:P\rightarrow G$ be the inclusion of a parabolic subgroup. For
every $r\in \mathcal{F}_{G}\cap \mathcal{F}_{P}$, Lemma 3.2 provides a basic
Weyl invariant $I\in H^{2r}(B_{G};\mathbb{Q})$ such that

\begin{center}
\noindent $q_{P}(i^{\ast }(I))\neq 0$ in $Q^{2r}(P)$,
\end{center}

\noindent where $i^{\ast }:H^{\ast }(B_{G};\mathbb{Q})\rightarrow H^{\ast
}(B_{P};\mathbb{Q})$ is the map induced by $i$. Since $q_{G}(I)\in Q^{2r}(G) 
$ and

\begin{center}
$Q^{2r}(i)(q_{G}(I))=q_{P}(i^{\ast }(I))$,
\end{center}

\noindent the induced map

\begin{center}
$Q^{2r}(i):$ $Q^{2r}(G)\rightarrow Q^{2r}(P)$
\end{center}

\noindent is nonzero by Lemma 3.2. Taking $s=2r$ in the above naturality
diagram, we conclude that

\begin{center}
$i_{\ast }\otimes 1$ $:\pi _{2r-1}(P)\otimes \mathbb{Q\rightarrow }\pi
_{2r-1}(G)\otimes \mathbb{Q}$
\end{center}

\noindent is nonzero. This proves the theorem.$\square $

\bigskip

\noindent \textbf{Remark 3.5.} Lemma 3.2 strengthens certain earlier results
concerning the rational homotopy groups of homogeneous spaces. In \cite{Me},
Meier reduced the study of the rational homotopy type of flag manifolds to
the construction of Sullivan models, without making the differentials
explicit. Our approach, by contrast, gives an explicit description of the
rational homotopy groups in terms of the rational homotopy types of $G$ and $%
P$, without requiring an explicit computation of these differentials.

In \cite[Proposition 3.9]{ST}, Shiga and Tezuka established a regularity
property for the restrictions of basic Weyl invariants. In the present
notation, their result implies the nonvanishing of the restriction of a
basic invariant of degree $2r$, for $r\in \mathcal{F}_{G}\cap \mathcal{F}%
_{P} $. Lemma 3.2 establishes a stronger conclusion: such an invariant $I$
can be chosen so that

\begin{center}
$q_{P}(I\mid _{P})\neq 0$.
\end{center}

\noindent In other words, the restriction remains nonzero in the
indecomposable quotient. This stronger nonvanishing property is essential to
the proof of Theorem 1.1 in Section 4.

More generally, Terzi\'{c} \cite[Proposition 13]{T} observed that, for a
common exponent of a compact connected semisimple group $G$ and a closed
connected subgroup $H$, the determination of the corresponding rational
homotopy groups depends on whether the relevant Cartan differential $%
d(z_{r}) $ vanishes. In particular, when the common exponent has
multiplicity one, the cases $d(z_{r})=0$ and $d(z_{r})\neq 0$ lead to
different rational homotopy groups. For the flag manifolds $G/P$ considered
here, Lemma 3.2 provides the additional nonvanishing information needed to
determine the relevant differential and hence to derive the explicit formula
in Theorem 1.1.$\square $

\section{Proof of Theorem 1.1}

Let $G$ be a simple Lie group. For every $r\in \mathcal{F}_{G}$ Table 1 gives

\begin{center}
$\pi _{2r-1}(G)\otimes \mathbb{Q=}\left\{ 
\begin{tabular}{l}
$\mathbb{Q\oplus Q}$ if $G=D_{n}\text{, }n\text{ }$even, $r=n$; \\ 
$\mathbb{Q}$ otherwise.%
\end{tabular}%
\right. $
\end{center}

\noindent Moreover, if $n\text{ is }$even and $n\in \mathcal{F}_{D_{n}}\cap 
\mathcal{F}_{P}$, then Table 3 implies that $P$ has local type with
semisimple part of type $A_{n-1}$ or $D_{m}$, where

\begin{center}
$\frac{1}{2}+1\leq m\leq n-1$.
\end{center}

\noindent Consequently, Theorem 3.4 admits the following refinement.

\bigskip

\noindent \textbf{Corollary 4.1. }For every $r\in \mathcal{F}_{G}\cap 
\mathcal{F}_{P}$, the induced map

\begin{center}
$i_{\ast }:\pi _{2r-1}(P)\otimes \mathbb{Q\rightarrow }\pi _{2r-1}(G)\otimes 
\mathbb{Q}$
\end{center}

\noindent has the following properties:

(a) If $(G,r)\neq (D_{n},n)$ with $n$ even, then

\begin{center}
$\pi _{2r-1}(G)\otimes \mathbb{Q\cong Q}$
\end{center}

\noindent and $i_{\ast }$ is surjective;

(b) If $(G,r)=(D_{n},n)$ with $n$ even, then

\begin{center}
$\pi _{2r-1}(P)\otimes \mathbb{Q\cong Q}$, $\pi _{2r-1}(G)\otimes \mathbb{%
Q\cong Q\oplus Q}$
\end{center}

\noindent and $i_{\ast }$ is injective.\noindent $\square $

\bigskip

We are now ready to prove Theorem 1.1, stated in Section 1.

\bigskip

\noindent \textbf{Proof of Theorem 1.1.} By (1.1) the rational homotopy
groups of $P$ and $G$ vanish in even degrees; that is

\begin{center}
$\pi _{2s}(P)\otimes \mathbb{Q}=\pi _{2s}(G)\otimes \mathbb{Q}=0$, $s\geq 1$.
\end{center}

\noindent Consequently, the rational homotopy long exact sequence associated
with the fibration

\begin{center}
$P\hookrightarrow G\rightarrow G/P$
\end{center}

\noindent breaks up into short exact sequences

\begin{enumerate}
\item[(4.1)] $0\rightarrow \pi _{2s}(G/P)\otimes \mathbb{Q\rightarrow }\pi
_{2s-1}(P)\otimes \mathbb{Q}\overset{i_{\ast }}{\mathbb{\rightarrow }}\pi
_{2s-1}(G)\otimes \mathbb{Q\rightarrow }\pi _{2s-1}(G/P)\otimes \mathbb{%
Q\rightarrow }0$.
\end{enumerate}

\noindent We determine the rational homotopy groups of $G/P$ by considering
the following cases.

\begin{enumerate}
\item[(i)] Suppose that $s\notin \mathcal{F}_{G}$ or $s\notin \mathcal{F}_{P}
$, then $i_{\ast }=0$, and (4.1) gives the isomorphisms
\end{enumerate}

\begin{center}
$\pi _{2s}(G/P)\otimes \mathbb{Q}\cong \pi _{2s-1}(P)\otimes \mathbb{Q}$, 

$\pi _{2s-1}(G)\otimes \mathbb{Q\cong }\pi _{2s-1}(G/P)\otimes \mathbb{Q}$.
\end{center}

\begin{enumerate}
\item[(ii)] Suppose next that $r\in \mathcal{F}_{G}\cap \mathcal{F}_{P}$. By
Corollary 4.1, we have two possibilities.
\end{enumerate}

(a) If $G=D_{n}$, $r=n$, and $n$ is even, then (4.1) becomes

\begin{center}
$0\rightarrow \pi _{2r}(G/P)\otimes \mathbb{Q}\rightarrow \mathbb{Q}\overset{%
i_{\ast }}{\mathbb{\rightarrow }}\mathbb{Q\oplus Q}\rightarrow \pi
_{2r-1}(G/P)\otimes \mathbb{Q}\rightarrow 0$,
\end{center}

where $i_{\ast }$ is injective. It follows that

\begin{center}
$\pi _{2r}(G/P)\otimes \mathbb{Q}=0$, $\pi _{2r-1}(G/P)\otimes \mathbb{%
Q\cong Q}$;
\end{center}

(b) In all other cases, (4.1) becomes

\begin{center}
$0\rightarrow \pi _{2r}(G/P)\otimes \mathbb{Q}\rightarrow \pi
_{2r-1}(P)\otimes \mathbb{Q}\overset{i_{\ast }}{\mathbb{\rightarrow }}%
\mathbb{Q}\rightarrow \pi _{2r-1}(G/P)\otimes \mathbb{Q}\rightarrow 0$
\end{center}

where $i_{\ast }$ is surjective. Hence

\begin{center}
$\dim _{\mathbb{Q}}\pi _{2r}(G/P)\otimes \mathbb{Q}=\dim _{\mathbb{Q}}\pi
_{2r-1}(P)\otimes \mathbb{Q}-1$, 

$\pi _{2r-1}(G/P)\otimes \mathbb{Q}=0$.
\end{center}

In particular, in both cases (a) and (b), we have

\begin{center}
$\dim _{\mathbb{Q}}\pi _{2r}(G/P)\otimes \mathbb{Q}=\dim _{\mathbb{Q}}\pi
_{2r-1}(P)\otimes \mathbb{Q}-1$.
\end{center}

\noindent Combining the above cases with convention (1.2), we obtain formula
(1.3), as asserted.$\square $

\bigskip

Haibao Duan, dhb@math.ac.cn

Yau Mathematical Science Center, Tsinghua University, Beijing 100084;

Academy of Mathematics and Systems Sciences, Chinese Academy of Sciences,
Beijing 100190.

\end{document}